\documentclass{amsart}
\usepackage[utf8]{inputenc}
\usepackage{enumerate}
\usepackage{amsmath, amsthm, amssymb, latexsym, graphics}
\usepackage{parskip}

\def \bsigma{{\boldsymbol{\sigma}}}
\def \bvfi{{\boldsymbol{\varphi}}}

\def \Fix{\operatorname{Fix}}

\newcommand{\tq}{\;\: | \:\:}
\newcommand{\mesp}{\;\:}

\newtheorem{theorem}{Theorem}[section]
\newtheorem{lemma}[theorem]{Lemma}

\newtheorem{proposition}[theorem]{Proposition}
\newtheorem{corollary}[theorem]{Corollary}
\newtheorem*{intth}{Theorem}

\newtheorem{definition}[theorem]{Definition}

\newtheorem{remark}[theorem]{Remark}
\newtheorem{fact}[theorem]{Fact}

\theoremstyle{definition}
\newtheorem{example}[theorem]{Example}

\AtBeginDocument{%
  
  {}}

\AtEndDocument{\bigskip{\footnotesize%
  \textsc{ Instit\"{u}t f\"{u}r Mathematishe Logik und
Grundladenforchung Fachbereich und Informmatik 
  Einsteinstra{\ss}e 62, 48149 M\"{u}nster, Deutschland
} \par 
  \textit{E-mail address}, G.Onay: \texttt{onay@uni-muenster.de} }}
\renewcommand{\epsilon}{\varepsilon}

\renewcommand{\leq}{\leqslant}
\renewcommand{\geq}{\geqslant} 

\newcommand{\jump}[1]{\operatorname{Jump}(#1)}
\newcommand{\jumpm}[2]{\operatorname{Jump}_{#1}(#2)}

\newcommand {\vfi}{\varphi}

\newtheorem{prop}[theorem]{Proposition}

\newtheorem{corr}[theorem]{Corollary}

\newtheorem{defn}[theorem]{Definition}

\theoremstyle{definition}
\newtheorem{exemple}[theorem]{Example}

\newtheorem*{nota}{Notation}

\usepackage{hyperref}
\hypersetup{
hidelinks=true}

\usepackage[left=3cm,right=3cm,top=3cm,bottom=3cm]{geometry}

\author{Gönenç Onay}
 \thanks{Research partially funded by DFG through SFB 878.}

\makeatletter
\def\blfootnote{\gdef\@thefnmark{}\@footnotetext}
\makeatother
\title{Valued Modules over Skew Polynomial Rings 2}

\begin{document}
\begin{abstract}
	Following our first article, we continue to investigate ultrametric modules  over  a  ring of twisted polynomials of the form $[K;\vfi]$, where $\vfi$ is a ring endomorphism of the field $K$. The main motivation comes from the the theory of valued difference fields (including characteristic $p>0$ valued  fields  equipped with the Frobenius endomorphism). We introduce the class of modules, that we call, affinely maximal and residually divisible and  we prove (relative -) quantifier elimination results. Ax-Kochen \& Erhov type theorems follow. As an application, we axiomatize, as a valued module,  the ultraproduct of algebraically closed valued fields $(\mathbb{F}_{p^n}(t)^{alg})_{n\in \mathbb{N}}$, of fixed characteristic $p>0$,  each equipped with the morphism $x\mapsto x^{p^n}$ and with the $t$-adic valuation.
\end{abstract}

\maketitle

{\small {\it MSC\ }: Primary 03C60; Secondary 03C10, 14G17, 12J20, 16W80, 13J15
\medskip

{\it Keywords\ }: Valued modules, Valued difference fields, Valued fields of positive characteristic, Quantifier Elimination, Ax-Kochen and Ershov Theorems.}

\section{Introduction}

Following our first article \cite{Onay2017}, we extend our definition of valued  modules over non-commutative rings of twisted polynomials and study their model theory.

Model theory of valued fields (including additional structures, such as endowed with an automorphism or a derivation) is much less understood in positive characteristic then it is in characteristic  $0$. For instance, an analog of famous Ax-Kochen and Ershov Theorem is not known for henselian valued fields of equicharacteristic $(p,p)$ $(p>0)$; for a survey see \cite{kuhlmanntame} and see our result  \cite{Onay2018} for the module structure of $\mathbb{F}_p((X))$. Notice also that, the theory of non-standard Frobenius in the valued setting, that is the theory of the field	$$\mathbf{F}_0:=\prod_{p \ \leadsto \ U} \left(\mathbb{F}_{p}(t)^{alg}, v_{\mbox{t}-adic}, x\mapsto x^{p} \right)$$ equipped with the {\it non-standard} Frobenius  $\vfi:=\lim_{p \leadsto U} x\mapsto x^{p}$ (where $U$ is a non principal ultrafilter over  prime numbers) is axiomatized by \cite{azgin-contractive} (as a valued difference field) but the theory of its positive characteristic analog given in the abstract is still unknown (as a valued difference field).

Here we propose a systematic study of -in particular-  the additive reducts of valued difference fields.

 Recall that, a valued difference field $(L,\vfi,v)$ is a field  
$L$ equipped with a valuation $v$, and an endomorphism $\vfi$ of $L$ such that $\vfi ({\mathcal{O}}) \subseteq \mathcal {O}$, where $\mathcal{O}$ is the valuation ring of $v$. Note that 
the theory of difference fields in fixed characteristic, $p=0$ or $p$ is a prime, admits a model companion denoted as $\operatorname{ACFA}_p$ by \cite{hu-frob} and the theory of valued difference fields of residue characteristic $0$ admits a model companion by \cite{azgin-contractive}.

 Let $K$ be a difference subfield of $(L,\vfi,v)$. We aim to study the $K$-vector space structure of $L$ together with the induced structure by $\vfi$ and by $v$. Here the $K$-vector space structure of $L$ together with $\vfi$ can be summarized as the right $R$-module structure of $L$ where
$R=K[t, \vfi]$, is the skew polynomial algebra over $(K,+)$, (right-) twisted by $\vfi$ respecting  the commutation rule $at=ta^{\vfi}$ for all $a\in K$. Interpreting $x.t$ as 
$x^{\vfi}$ in $K$, makes $K$, a $K[t;\vfi]$-module. Note that every $r\in R$ can be expressed in the form $r=\sum_{i} t^ia_i$, and this expression is unique (for more details see \cite{Onay2017}, section 2).

For a $R$ as above, we consider two-sorted structures 
$(M, v:M\twoheadrightarrow \Delta)$, where $M$ is the (right-) $R$-module sort, $\Delta$ is the linear order sort equipped with a top element $\infty$  and with an (right-) action of $R$, denoted as $\delta\cdot r$, for $\delta\in \Delta$ and $r\in R$. We call $\Delta$, equipped with the action of $R$, an $R$-chain. Moreover, we require that
$(M,+,v)$ be an valued abelian group such that $v(x.ta)=v(x)\cdot ta$ for all $x\in M$ and $a\in K$. More generally, scalar multiplication $M\times R \to M$, is required to be compatible with the $R$-chain structure on $\Delta$, on {\it generic} points: for instance, for a fixed $r \in R$, 
\begin{equation}
v(x.r)=v(x)\cdot r \tag{*} 
\end{equation} for all $x$ such that $v(x)$ is not in a finite set of exceptional values, depending only on the action of $r$  on $\Delta$. These exceptional values called {\it potential jumps} of $r$, are analogs of tropical zeros of a usual polynomial.

A point $x$ satisfying ({*}) is called a {\it regular} for $r$. Consider the case where $K$ is an algebraically closed field of characteristic $p>0$, with $\vfi:x\mapsto x^p$ and $R=K[t;\vfi]$. The ring $R$ can be seen as the ring of additive polynomials equipped with the composition and the addition. Then, $x\in F\subseteq K$,  is regular for all non zero $r \in R$, if and only if the type of $x$ over $K$, in the language of valued fields, is strongly stably dominated (see \cite{hu-loeser}, Chapter 8). Regularity is one of the central notions of the present article.

We first establish a theory of pseudo-convergence \`a la Kaplansky and then introduce the class that we call {\it affinely maximal  $R$-modules}, an analog of algebraically maximal Kaplansky fields. This class is characterized with the simple fact that for every non zero $x$, for every non zero $r\in R$, there is some $y$, {\bf regular} for $r$ such that $y.r=x$. 

Let $L_{Mod}(R)$ be the language of (right-) $R$-modules and $L_V$ be the language of $R$-chains, that is the language of ordered sets together with a function symbol for each $r\in R$ and with a constant symbol $\infty$.  Set 
$$L:=L_{Mod}(R)\cup L_V\cup\{v\}.$$ Our Theorem \ref{thm:main} establishes a quantifier elimination result modulo the $\Delta$-sort from which we deduce the following results. 

\begin{theorem}[A-K,E $ \equiv $]\label{akeequiv}
	Let $ (F, v) $ and $ (G, w) $ be two non zero affinely maximal residually divisible modules, such that $F_{tor}$ and $G_{tor}$ are elementary equivalent as $L_{Mod}(R)$-structures and $ v (F) $ and $ v (G) $ are elementary equivalent in the language $ L_{V} $. Then $ (F, v) $ and $ (G, v) $ are elementary equivalent as $ L $-structures.
\end{theorem}

\begin{theorem}\label{thm:acfa}
	Let $(K,v_K, \vfi)$ be a valued difference field such that $(K,\vfi) \models \operatorname{ACFA}_{p}$ and $v_K(x^{\vfi})>nv(x)$ for all $x$ such that $v_K(x)>0$. The theory of affinely maximal residually divisible $K[t;\vfi]$-modules together the theory of dense $R$-chains is complete, eliminates quantifiers in $L$, and the valued $R$-module $(K,v_K)$ is the prime model of this theory.  
\end{theorem} 

Note that our results are optimal in the following sens. First, as we have already noticed in our first article, divisibility alone, while analog of algebraic closeness, does not imply good valuation theoretical properties in the case of modules (see \cite{Onay2017}, Proposition 3.26). Moreover some strong results known for valued fields, can not be obtained in our case (see \cite{gonenc}, Exemple 4.3.22).

This manuscript is a short but richer version (in term of results) of Chapter 4 in \cite{gonenc}, and is organized as follows. In section \ref{R-chains} we define the action of $R$ on linear orders and axiomatize the $L_V$-theory that we call the theory
full $R$-chains, and show that together with density axioms, this theory has quantifier elimination. In section \ref{v-mod}, we introduce valued $R$-modules, and axiomatize affinely maximal residually divisible modules after establishing a Kaplansky theory of pseudo-convergence. Some applications to the theory of valued fields in positive characteristic is given (see for instance \ref{corr:appfields}). In section 4, we establish embedding theorems and deduce the corresponding quantifier elimination results.

In order to keep our article in a suitable size, we have omitted some proofs that are already detailed in \cite{gonenc}. 

 Note also that we are preparing another article which deals with the geometric consequences of our results in terms of minimality notions in model theory (in particular C-minimality).

{\bf Acknowledgments.} I'm grateful to Françoise Delon who insistently encouraged me to publish this manuscript. 

\section{R-Chains}\label{R-chains}

During this section, we let the pair $(K,\vfi)$ range over  difference fields, that is  $K$ is a field together with an ring endomorphism $\vfi$ and we let
$R:=K[t;\vfi]$.

\begin{defn}[$K$-chains] For $a\in K$, let $\cdot a$ denotes a unary function symbol (acting on the right).
  A $K$-chain, is a structure of the form $(\Delta,<,\infty, (\cdot a)_{a\in K})$ where $(\Delta, <, \infty)$ is a linear order with a top element $\infty$ and such that  
  for all  non zero $a,b \in K$ and $\gamma, \delta \in
  \Delta \setminus \{\infty\}$, we have:
  \begin{enumerate}
  \item $\gamma > \delta \to \gamma \cdot a> \delta \cdot a$,
  \item $\gamma \cdot ab = (\gamma \cdot a) \cdot b$,
  \item $\gamma \cdot a > \gamma \to \delta 
  \cdot a > \delta$,
  \item $\gamma \cdot (a \pm b)\geq 
  \min\{\gamma\cdot a, \gamma \cdot b\}$,
  \item $\gamma \cdot 0= \infty \; ;
  \gamma \cdot 1 =\gamma\; ; \infty\cdot b=\infty$.
  \end{enumerate}
\end{defn}
\begin{remark}
	If $a \in K$, and $\gamma$ is not infinity, then 
	$\gamma \cdot a=\infty$ implies  $a=0$.
\end{remark}

\begin{nota}
	In the rest of this paper, whenever $(M,v)$ is a valued structure, a valued field, valued module, valued abelian group etc., for a subset $A\subseteq M$ we denote by $vA$, the set
	$\{v(a)\tq a\in A\}$. 
\end{nota}
\begin{example} If $(K,v)$ is a valued field then
	letting, for $\lambda\in K$, $v(\lambda) \cdot a := v(\lambda a)= v(\lambda)+v(a)$ makes the ordered set 
	$vK$ a $K$-chain.
\end{example}

 Conversely we have the following.

\begin{proposition}\label{valinK}
	If $\Delta$ is an infinite $K$-chain then 
	$$\mathcal{O}_{\Delta}:=\{a \in K \tq \gamma \cdot a \geq \gamma \; \text{for all} \; \gamma\neq \infty \},$$ 
	is a valuation ring of $K$.
\end{proposition}
\begin{proof}
	Same proof as the one given in \cite{fares}, Proposition 10.
\end{proof}
\begin{nota} When the $K$-chain $\Delta$ is clear from the context,  we denote by $v_K$ the valuation on  $K$ given by    $\mathcal{O}_\Delta$. 
\end{nota}

\begin{theorem}\label{fares-elimination} Set the language $\mathcal{L}:=\{<,\cdot a \ (a\in K), c, \infty\}$. Then  theory of dense $K$-chains  such that $\Delta\setminus \{\infty\}$ is  without end points, together with the axiom $\infty \neq c$ is complete and eliminates quantifiers in the language
	$\mathcal{L}$.  
\end{theorem}
\begin{proof}
	Follows immediately by Th\'eor\`eme 6, and Th\'eor\`eme 13 in \cite{fares}. 
\end{proof}

\begin{nota}For $r \in R$, let $\cdot r$ be a unary function symbol.
	We set $L_V:=\{<,\cdot r_{r\in R}, \infty\}$, that we call the language of $R$-chains.  
\end{nota}

Recall that any $r\in R$ can be uniquely written as $\sum_{i} t^ia_i$, with $a_i\neq 0$. 
We call a term $t^ia_i$ in this expression {\bf a monomial} of $r$.

\begin{defn}[$R$-chains]\label{rchains}
  An $R$-chain  $(\Delta, <, \infty, \cdot r_{r \in R})$ is an a $K$-chain such that
  \begin{enumerate}
  \item $\cdot t$ is strictly increasing on $\Delta \setminus \{\infty\}$, and $\infty \cdot t=\infty$,
 \item $\gamma\cdot ta=(\gamma \cdot t)\cdot a$ and $\gamma\cdot t^na = (\gamma \cdot t^{n-1})\cdot ta$ for all $\gamma \in \Delta$, and $a\in K$,
  \item $\gamma \cdot r = \min_{i} \{\gamma \cdot \mathbf{m}_i\}$, for all non zero $r \in R$, where
  the $\mathbf{m}_i$ are the monomials of $r$ and for all $\gamma \in \Delta$,
  \item $\gamma \cdot \mathbf{m}_1 \leq \gamma \cdot \mathbf{m}_2 \to \left( 
  \forall \delta \mesp (\delta <\gamma \to \delta \cdot 
  \mathbf{m}_1 < \delta \cdot \mathbf{m}_2)\right)$ for all monomials 
  $\mathbf{m}_1, \mathbf{m}_2 \in R$ such that $0\leq\deg(\mathbf{m}_2)<\deg(\mathbf{m}_1)$,
		  and for all $\gamma
		  \neq \infty$.
  \end{enumerate}
\end{defn}

Note that by above axioms for every non zero monomial $\mathbf{m}$, 
$\gamma \mapsto \gamma \cdot \mathbf{m}$ is strictly increasing and the axiom scheme $4.$ implies
its dual:
\begin{equation}
\gamma \cdot \mathbf{m}_1 \leq \gamma \cdot \mathbf{m}_2 \to \left( 
  \forall \delta \mesp (\delta >\gamma\to \delta \cdot 
  \mathbf{m}_1 < \delta \cdot \mathbf{m}_2)\right) \tag{$^d4$}
\end{equation}
 for all monomials 
 $\mathbf{m}_1, \mathbf{m}_2$ such that $0\leq\deg(\mathbf{m}_1)<\deg(\mathbf{m}_2)$.

\begin{remark}\label{tausurjective} \begin{enumerate}
	\item If $\cdot t$ is onto, then for every non zero $\delta \in \Delta$, and for every non zero $r\in R$, 
	there is a (unique) $\gamma \in \Delta$, such that $\gamma \cdot r =\delta$. 
	
	\item We will denote by 
	$\cdot t^{-k}$ the inverse of $\cdot t^k$ for $k\in \mathbb{N}_{>0}$. Suppose also that $\vfi$ is an automorphism  of $K$. Then the inequation 
	$$\gamma \cdot t^ja\leq \gamma \cdot t^i b\mesp (ab\neq 0),$$ where $j\geq i$, is equivalent to
	$$\gamma \cdot t^{j-i}\vfi^{-n}(a/b)\leq \gamma.$$
	
\end{enumerate} 
\end{remark}

Recall that a valued difference field is a valued field $(K,v)$ together with an field endomorphism $\vfi$, such that 
$\vfi(\mathcal{O})\subseteq \mathcal{O}$, where $\mathcal{O}$ is the valuation ring of $v$. It follows that 
$\vfi$ induces endomorphisms $\bvfi$ and $\bar{\vfi}$, respectively on the value group and on the residue field. A difference polynomial map on $K$, is of the 
form $x \mapsto f(x,\vfi(x),\dots,\vfi^n(x))$ where 
$f(x_0, \dots, x_n)$ is an ordinary polynomial  over $K$ and we denote it as $f^{\vfi}:K\to K, \mesp x \mapsto f^{\vfi}(x)$.
  For every difference polynomial map $f^{\vfi}$, we associate the {\it tropical} map $f^{\vfi}_v:vK\to vK$, given by 
$$\gamma \mapsto Trop(f)(\gamma,\bvfi(\gamma),\ldots, \bvfi^n(\gamma)),$$
where $Trop(f)(\gamma_0, \ldots, \gamma_n)=\min_i\{v(\mathbf{m}_i(a_0,\ldots,a_n))\}$, the
$\mathbf{m}_i$ range over monomials of $f$, and $$(v(a_0),\ldots, v(a_n))=(\gamma_0, \ldots, \gamma_n).$$ Note that  $Trop(f)(\gamma_0, \ldots, \gamma_n)$ depends only on  the tuple of valuations  $(v(a_0),\ldots, v(a_n) )$. 
We say that $f^{\vfi}$ is linear when $f$ is. When $\vfi$ is not of finite order in $\operatorname{End}(K)$, the ring  $R:=K[t;\vfi]$ is isomorphic to ring of linear difference maps equipped with the composition and the addition (see \cite{cohn} and \cite{Durhan2015} for more details). 

{\bf Assumption/Notation.} For the rest of this article we consider only difference fields $(K,\vfi)$, such that $\vfi$ is not of finite order in $\operatorname{End}(K)$. For $r\in R$, we write $r_v$, to denote the tropical map associated to the linear difference  polynomial $r^{\vfi}$ which is the image of $r$ in $\operatorname{End}(K)$.

\begin{exemple} Let $(L,\sigma, v)$ be a valued difference field extension of $(K, \vfi, v_K)$. Denote by $\bsigma$ the endomorphism induced by $\sigma$ on 
$vL$. Assume that $\bsigma(\gamma)>\gamma$ for $\gamma > 0$ (this assumption is needed for the axioms (4)). For $r=t^na_n+\ldots+a_0 \in R$ set
the interpretation $$ \gamma\cdot r := r_v(\gamma)=\min_{i=1,\dots,n} \{\bsigma^i(\gamma)+v(a_i)\}.$$
Then $vL$ is an $R$-chain. 
\end{exemple}

Conversely we have:
\begin{remark} Let $\Delta$ be an infinite $R$-chain, such that $\gamma \mapsto \gamma \cdot t$ is surjective on $\Delta$. Let $v_K$ be the valuation on $K$, induced by the $K$-chain structure of $\Delta$. Then $(K,v_K,\vfi)$ is a valued difference field.  
\end{remark}
\begin{proof}
Let $\mathcal{O}$ be the valuation ring of $v_K$ and let $\lambda \in \mathcal{O}$.  Let $\delta$
be such that $\delta\neq \infty$, and take $\gamma$ such that $\delta=\gamma \cdot t$. Since $\lambda  \in \mathcal{O}$, and $\cdot t$ is strictly increasing, 
$\delta=\gamma\cdot t\leq (\gamma \cdot \lambda)\cdot  t =\gamma\cdot t\lambda^{\vfi}=\delta\cdot \lambda^{\vfi}$. That is, $\lambda^{\vfi} \in \mathcal{O}$.
\end{proof} 

Let $\Delta$ be an non-empty $R$-chain. The following lemmas are easy consequences of Definition \ref{rchains}.

\begin{lemma}\label{lem:consq:facile}
For all $\gamma \in \Delta$ and $r,q \in R$ we have, 
$\gamma \cdot (r\pm q) \geq \min\{\gamma \cdot r, \gamma \cdot q\}$ and 
$\gamma \cdot r \neq \gamma \cdot q$ then 
$\gamma \cdot (r\pm q)=\min\{\gamma \cdot r, \gamma \cdot q\}$.
\end{lemma}
\begin{proof}
	This is straightforward.
\end{proof}

\begin{lemma}\label{lem:consq:cuts}
Let $\mathbf{m}_1$ and $\mathbf{m}_2$ be non zero monomials from $R$, such that $0\leq \deg(\mathbf{m}_2)<\deg(\mathbf{m}_1)$. Then there exists at most one $\gamma\neq \infty$ such that 
$\gamma \cdot \mathbf{m}_1 =\gamma \cdot \mathbf{m}_2$ and the set 
$$A:=\{\gamma \in \Delta \tq \gamma\cdot \mathbf{m}_1 < \gamma \cdot \mathbf{m}_2\}$$ is an initial segment 
 of $\Delta$. Suppose
 $\cdot t $ is onto, $A$ is proper non-empty subset of $\Delta\setminus \{\infty\}$ and there exists no $\gamma\neq \infty$ such that 
 $\gamma \cdot \mathbf{m}_1 =\gamma \cdot \mathbf{m}_2$.  Then $A$ defines a Dedekind cut; if such a $\gamma$ exists then $\gamma  =\sup A = \inf ( (\Delta\setminus\{\infty\}) \setminus A )$.  
\end{lemma}
\begin{proof}
	The fact that $A$ is an initial segment (eventually empty), and the uniqueness of a such $\gamma$  follows directly from Axioms (4) of the definition \ref{rchains}.
	 
	Suppose $\cdot t$ is onto, and $A\neq \emptyset$. 
	Let $\delta \in A$. Write $\mathbf{m}_1=t^ja$ and
	$\mathbf{m}_2=t^ib$ with $j>i$. For sure we may suppose $a=1$. Let $\delta \in A$,  then $\delta < (\delta\cdot t^ib )\cdot t^{-j} \in A$. Hence $A$ has no maximum. Similarly one shows that the complementary of $A$ has no minimum. 
\end{proof}

\begin{defn}[Potential jump values]\label{def:potetialjumps}
	For $r \in R$ we call a potential jump of $r$ in $\Delta$,  any element $\gamma\neq \infty$ such
	that $$\gamma \cdot \mathbf{m}_i= \gamma \cdot \mathbf{m}_j=\gamma \cdot r$$ for different monomials  $\mathbf{m}_i$ and $\mathbf{m}_j$ of $r$ (in particular  $\mathbf{m}_i$ and $\mathbf{m}_j$ have different degrees).
	We denote by $\jumpm{\Delta}{r}$ the set of potential jump values of $r$
	in $\Delta$.
	
	A potential jump of an $R$-chain $\Delta$, 
	is a potential jump for some non zero $r$. We denote by 
	$\jumpm{\Delta}{R}$ the set of all potential jumps in $\Delta$.
\end{defn}
\begin{remark}\label{rem:potjumps}
	\begin{enumerate}
\item	$\jumpm{\Delta}{R}$ is an $L_V$-substructure of $\Delta$.
\item 	Potential jumps of in $\Delta$,  can also be defined as the set of $\infty\neq \gamma \in \Delta$ such that  
	$$\gamma \cdot \mathbf{m_1} = \gamma \cdot \mathbf{m}_2$$ for
	all monomials $\mathbf{m}_1,\mathbf{m}_2$ satisfying $0\leq\deg(\mathbf{m}_1)<\deg(\mathbf{m}_2)$.

\item\label{dcl-r-chain} Let $\operatorname{dcl}^{\Delta}$ denotes the definable closure operator in an $R$-chain $\Delta$. Then $\jumpm{\Delta}{R} \subseteq \operatorname{dcl}^{\Delta}(\emptyset)$.
\end{enumerate}
\end{remark}

\begin{corollary}\label{cor:jumpsiso}
	Let $\Delta$, $\Delta'$ be $R$-chains such that
	$\jumpm{\Delta}{R} \equiv \jumpm{\Delta'}{R}$ as $L_V$-structures. Then $\jumpm{\Delta}{R}$ and 
	$\jumpm{\Delta'}{R}$ are isomorphic. 
\end{corollary}
\begin{proof}
	Follows by the 3th claim above.
\end{proof}

Recall for the following definition that $R$ embeds in a smallest division ring, denoted by $D$ as in our first article (see \cite{Onay2017}). Suppose $\cdot t$ is onto, then by Remark
\ref{tausurjective} (1), the action of $R$ on $\Delta$, extends to an action of $D$ on $\Delta$.

\begin{defn}\label{def:fullchains} 
  	$\Delta$ is said to be full if
  \begin{enumerate}
  \item $\gamma \mapsto \gamma \cdot t$ is onto,
  \item every subset of the form $\{\gamma \in \Delta \tq \gamma \cdot \mathbf{m} \ \bullet \
  \gamma, \gamma\neq \infty\}$ of $\Delta$, where $\mathbf{m} \notin K$ a monomial and $\bullet\in \{<,>\}$, is a non empty proper subset of $\Delta \setminus \{\infty\}$,
  \item $\Delta\setminus\{\infty\}$ contains a solution for all equation of the form 
  $\gamma\cdot t^n=\gamma \cdot a$ for all non zero $a\in K$ and $n\in \mathbb{N}_{>0}$.
  \end{enumerate}
  We say that $\Delta$ is a $D$-chain, if $\Delta$ satisfies (1) and (2). 
\end{defn}
\begin{remark}\label{rem:inducesR-chainonvK}
	Suppose $\Delta$ is full. Let $\theta\in \Delta$ be the unique
	element such that $\theta\cdot t =\theta$ (uniqueness follows by Lemma \ref{lem:consq:cuts}). Then $
	v_KK$ embeds in $\Delta$ via the map $v(\lambda) \mapsto \theta\cdot \lambda$ and inherits a structure of $R$-chain by letting $v(\lambda)\cdot r := \theta \cdot \lambda r$. In particular, $v(\lambda)\cdot t=(\theta\cdot \lambda t) =\theta\cdot(t\lambda^{\vfi})=(\theta \cdot t)\cdot \lambda^{\vfi}=\theta\cdot \lambda^{\vfi}$. It follows that
	$v(\lambda)\cdot r=\theta \cdot r_v(v(\lambda))$. 
\end{remark}

{\bf Assumption.} For the rest of the paper we assume that for any $R$-chain (in particular any $D$-chain) that we consider, the induced valuation $v_K$ on $K$ is non trivial.
\begin{corollary}[of Lemma \ref{lem:consq:cuts}]\label{densityofR-chains-1} For any $D$-chain  $\Delta$,   $\jumpm{\Delta}{R}$  is a densely ordered subset of $\Delta$. 
\end{corollary}
\begin{proof}
	Suppose $\gamma \in \jumpm{\Delta}{R}$. Write the equation $\gamma\cdot t^j=\gamma\cdot t^ib$ (with $j>i$). Let $\delta \in \jumpm{\Delta}{R}$ such that $\delta<\gamma$. By the proof of \ref{lem:consq:cuts}, $\delta<(\delta  \cdot t^{i}b)\cdot t^{-j}<\gamma$. Since $\jumpm{\Delta}{R}$ is a substructure of $\Delta$, 
	$\delta \cdot t^ib$ is potential jump of $\Delta$. It is clear then that $\delta \cdot t^ib\cdot t^{-j}$, is a also a potential jump.
\end{proof}

\begin{nota}Let $r=\sum \mathbf{m}_i$ where $\mathbf{m}_i$ are the monomials of $r$. For each $i$ we set 
	$$U_i(r):=\{\gamma \in \Delta \tq \gamma \cdot r = \gamma \cdot \mathbf{m}_i \}.$$ 
	Moreover we set 
	$J(r):=\{i \tq U_i \neq 0\}$
\end{nota}
\begin{remark}
	$J(r)$ is never empty and the union of $U_j(r)$ $(j\in J(r))$ is the whole chain $\Delta$.
\end{remark}

Fix a non zero $r$ and write $J(r)$ as $\{j_0,\ldots, j_k\}$ with 
$$m=j_0 < \ldots <j_k=n.$$
\begin{lemma}\label{lem:U_i}
	We have 
	$$ U_n(r)\leq \ldots U_{j_{i+1}} \leq U_{j_{i}} \ldots \leq U_m$$ where
	$U_n$ is an initial segment, $U_m$ is a final segment and each intermediate 
	$j_i$ is an bounded interval. Moreover, when
	$U_j\cap U_{j'}\neq \emptyset$, ($j\neq j'$), it is reduced to a singleton and for each $i$ such that $j<j_i<j'$,
	$U_{j_i}$ is equal to this intersection.
\end{lemma}
\begin{proof}
	See Lemme 4.1.12, in \cite{gonenc}.
\end{proof}
\begin{remark} Let $r=\sum_{i\in I} t^ia_i$. Set 
	$J(r)=\{j_0,j_1, \dots, j_k\}$ with $j_i<j_{i+1}$. Then
	$\jump{r}=\bigcup_{i=0..k-1} (U_{j_{i+1}}\cap U_{j_i})$.
\end{remark}
\begin{corollary}[of Lemma \ref{lem:U_i}]\label{cor:rchainmain}
	For all $r \in R \setminus\{0\}$, the map $\cdot r : \Delta \to \Delta$, $\gamma \mapsto \gamma\cdot r$ is strictly
	increasing. Moreover, for all $p, q \in R$ 
	we have  $(\gamma \cdot p) \cdot q = \gamma \cdot pq$.
\end{corollary}

\begin{proof}
	See Corollaire 4.1.14 in \cite{gonenc}.
\end{proof}

\begin{theorem}\label{ake-r-chains} Let $\Delta_0$ be a $D$-chain then $Th(\jumpm{\Delta_0}{R})$ eliminates the quantifiers in $L_V$.
\end{theorem}
	\begin{proof}
		Let $\Delta$ and $\Delta'$ both elementary equivalent to $\jumpm{\Delta_0}{R}$. Note in particular that $\Delta$ and $\Delta'$ are densely ordered and $\jumpm{\Delta}{R}\equiv
		\jumpm{\Delta'}{R}$ and hence by Corollary
		 \ref{cor:jumpsiso}, they are isomorphic.
		
		Let $\Gamma$ be a common substructure of $\Delta$ and $\Delta'$. By the above discussion we may suppose that $\Gamma$ contains all the potential jumps. Moreover it is clear that the $D$-chain generated by $\Gamma$ respectively in $\Delta$ and $\Delta'$ are isomorphic. Note that  $D$-chain generated by $\Gamma$ in $\Delta$, is the $R$-chain $$\{\delta \in \Delta, \tq \delta\cdot r \in \Gamma, \mesp \text{for some non zero}\mesp r\},$$ and  is a subset of $\operatorname{dcl}^{\Delta}(\Gamma)$. In short, we may assume that $\Gamma$ is a $D$-chain and contains all the potential jumps.

		Let $\phi(x,\bar{y})$ be a quantifier free $L_V$-formula and $\bar{a}\in \Gamma^{|\bar{y}|}$. We will show that $\phi(x,\bar{a})$ is satisfiable in $\Delta$ if and only if it is satisfiable in $\Delta'$. For this purpose,  we may assume that 
		$\phi$ is a disjunction of atomic or negation of atomic formulas. Write 
		$\bar{a}=(a_1, \ldots, a_n)$. Note that if $\phi(x,\bar{a})$ implies a formula of the form
		$x\cdot r=a_i\cdot q$, then, since $\Gamma$ is a $D$-chain, the unique element $\gamma \in \Delta$, verifying 
		$\gamma\cdot  r= a_i.q$ is already in $\Gamma$. Note also that any clause of the form 
		$x \cdot r =x \cdot q$ defines a potential jump, hence already in $\Gamma$. A clause of the form $x\cdot r < x\cdot q$, defines a non empty initial or final segment by the definition of $D$-chains, moreover  it is infinite without a least (respectively top) element, since $v_K$ is non trivial.

		We may hence assume that $\phi(x,\bar{a})$ is equivalent to a disjunction of formulas of the
			form $x\cdot r<a_i\cdot q_i$ or 
			$x\cdot q > a_j\cdot q_j$. Since $\Gamma$ is a $D$-chain, and $\gamma \mapsto \gamma \cdot r$ is strictly increasing by \ref{cor:rchainmain}, $\phi(x,\bar{a})$ defines a union intervals with end points in $\Gamma$.  Now, by density of both $\Delta$ and $\Delta'$, if $\phi(x,\bar{a})$ is satisfied by some $\delta$ in $\Delta$, then it is satisfied by some $\delta' \in \Delta'$. This finishes the proof. \end{proof}

\begin{corr}\label{D-chain-corr-eq}
$\jumpm{\Delta_0}{R}$ is the prime model of its theory.
\end{corr}

\begin{prop} Any $D$-chain embeds in
a full $R$-chain $\hat{\Delta}$ in which it is dense (with respect to the order topology), and  $\hat{\Delta}$ is up to an $R$-chain isomorphism.
\end{prop}
\begin{proof}
Set $\hat{\Delta}=\Delta\cup  CS$, where $CS$ is the set of Dedekind cuts defined in the proof of \ref{lem:consq:cuts} ordered by inclusion. It is straightforward to verify that is an $R$-chain and uniqueness follows by construction (see also Proposition 4.1.18 in \cite{gonenc}).
\end{proof}

\begin{definition}
We call $\hat{\Delta}$ as above the full closure of $\Delta$.
\end{definition}

\begin{corr}\label{cor:fullRchains}
The theory of dense and full $R$-chains is complete, eliminates quantifiers,  $o$-minimal and  is the model completion of the theory of dense $D$-chains.
\end{corr}

\begin{remark}\label{rem:acfa}
 Let $\Delta$ be full and dense $R$-chain. If $v_KK$ is already full, then $v_KK\preceq \Delta$. 
\end{remark}

\begin{example}\label{example-full-chains}
	\begin{enumerate}
		\item 	Any value group of an algebraically closed valued field of characteristic $p>0$, that is an ordered abelian divisible together with $\infty$, is a full and dense $K[t;x\mapsto x^p]$-chain (see below for the converse). 
		
		\item If $(K,v,\sigma)$ is a valued difference field, such that
		$v(\sigma(x))>v(x)$ for all $x$  with $v(x)>0$ and $vK$ is divisible as $\mathbb{Z}[\bsigma]$-module, then $vK$ is a full and dense $K[t;\sigma]$-chain. In particular this the case  when 
		$(K,\sigma) \models \operatorname{ACFA}_p$, and $\sigma$ is contractive: that is $v(\sigma(x))>nv(\sigma(x)$ for all $n\in \mathbb{N}$ and all $x$ such that $v(x)>0$. 
	\end{enumerate}

\end{example}

  \begin{remark} Let $p$ be a prime number and suppose that $(K,v)$ is of characteristic $p>0$. 
  Suppose $\Gamma:=vK \setminus \{\infty\}$ 
  equipped
  with the $K[t;x\mapsto x^p]$-chain structure is full, then $\Gamma$ is a divisible 
  abelian group.
  \end{remark}
  \begin{proof} Surjectivity of $\cdot t$ means that $\Gamma$ is $p$-divisible. Let $q\neq p$ be
  a prime number. Then for every $\delta$ there is $\gamma$ such that 
  $\gamma \cdot t^{q-1}=\gamma \cdot a $ where $v_K(a)=\delta$. In other words we have
  $p^{q-1}\gamma - \gamma = \delta$. Hence  $\Gamma$ is 
  $p^{q-1}-1$ divisible. By Little Fermat's Theorem, $\Gamma$ is $q$-divisible.
  \end{proof}

\section{Valued Modules}\label{v-mod}

We extend the notion of valued module given in \cite{Onay2017}. Recall that an valued abelian group $M$ is an abelian group together with an ordered set $\Delta$ having a top element $\infty$ and a surjective map
$v:M\to \Delta$ such that for all $x,y \in M$, 
$v(x\pm y)\geq \min\{v(x),v(y)\}$ and $v(x)=\infty$ if and only if $x=0$. 
\begin{defn} A valued $R$-module is a valued abelian group ($M, v:M\to \Delta$) where $\Delta$ is an $R$-chain
such that for all $x\in M$ one has
\begin{enumerate}
\item $v(x.a)=v(x)\cdot a$, for all $a \in K$,
\item  $v(x.t)=v(x)\cdot t.$
\end{enumerate}
\end{defn}

In the rest of the article, unless we specify explicitly the ring $R$, we say a valued module, instead of valued $R$-module. 

In a valued module $M$, we denote by $B(x,\gamma)$ 
the set $\{y\in M \tq v(x-y)>\gamma\}$ the open ball of radius $\gamma$ centered at $x$. Closed balls are defined similarly by changing $>$ to $\geq$. Recall that open balls from a clopen basis of a topology, called ultrametric topology.  

\begin{remark}
It follows that $v(x)\cdot \mathbf{m} = v(\mathbf{m}(x))$, for all monomial $\mathbf{m}$. Moreover, $v(x.r)\geq 
 \min_i\{v(x.\mathbf{m}_i)\}=\min_i \{v(x)\cdot \mathbf{m}_i\}=v(x)\cdot r$,
where $\mathbf{m}_i$ are the monomials of $r \in R$.
\end{remark}
Unless otherwise said, until the end of this section, we will denote by $(M,v)$ a valued $R$-module.

\begin{defn}[Regular Elements] For $r\in R$ and $x\in M$ we call $x$ regular for $r$ if 
$$v(x.r)=v(x)\cdot r .$$
 Otherwise $x$ is said to be   irregular for $r$.
If $x$ is regular for all $r$ then it is said to be regular for $R$. Moreover, 
if $A\subseteq M$ and $v(x.r-a)=\min\{v(x)\cdot r, v(a)\}$ for all $a$ and $r$ then 
$x$ is said to be regular over $A$. Finally if every $b\in B\subseteq M$ is regular (respectively regular over $A$) then $B$ is said to be
regular subset/subgroup/submodule (respectively regular subset/subgroup/submodule over $A$)
\end{defn}

\begin{remark} $0$ is regular for all $r$, and any $x$ is regular for any monomial and for the scalar $0_{R}$. 
\end{remark}
\begin{remark}\label{regularelementsexits}
Regularity is more generally investigated in \cite{Durhan2015}. Let $(K\subseteq M,v)$, be an extension of characteristic $p>0$ valued fields. Consider $(M,v)$ as valued module over $R:=K[t;x\mapsto x^p]$. It follows from Lemma 3.6 in \cite{Durhan2015} that,  whenever the residue field of $M$ is infinite, for every $r\in R$, for every $\gamma \in vM$, there is $x\in M$, regular for $r$ of valuation $\gamma$. 
\end{remark}

\begin{definition} The valuation $v(x)$ of an irregular $x\in M$ for $r\in R$, is called  a {\it jump value} of $r$ in $M$. 
\end{definition}
\begin{example}\label{examplejump}
	If $x$ non zero and $x.r=0$ for a non zero $r$, then $v(x)$ is a jump value of $r$ in $M$.
\end{example}
\begin{nota} We set 
$$\jumpm{M}{r}:=\{\text{jump values of}\mesp r \mesp \text{in}\mesp M\}.$$
\end{nota}

It follows easy from ultrametric triangle inequality that:
\begin{remark}
	$\jumpm{M}{r} \subseteq \jumpm{\Delta}{r}$.
\end{remark}

\begin{lemma}\label{rvregular} Let $x,y$ be such that $v(x-y)>v(x)$ then for any non zero
	$r$, is regular for $r$ if and only if, 
$y$ is.
\end{lemma}
\begin{proof} Note that $v(x)=v(y)$. Write the trivial equality $v(y.r)=v(x.r+(y-x).r)$. Since $\cdot r$ is strictly increasing 
$v(y)\cdot r = v(x) \cdot r < v(y-x)\cdot r\leq v((y-x).r)$. Hence $v(x.r)=v(x)\cdot r$ if and only if $v(y.r)=v(x.r)=v(y)\cdot r$. 
\end{proof}

The following is straightforward, one can also see \cite{gonenc}, Lemma 4.2.11.
\begin{lemma}
  For $r,q \in R$ and $x\in M$, we have
  \begin{enumerate}
  \item $x$ is regular for $r$ and $x.r$ is regular for $q$ if and only if 
  $x$ is regular for $rq$, 
  \item $(M,v)$ is non trivially valued and $|R|^{+}$-saturated then it contains
  non zero regular elements.
  \item The set of regular elements of $M$ is not stable under addition in general.
  \end{enumerate}
\end{lemma}

  \subsection{Extensions of valued modules}
  As the immediate extension of the valued  fields, we will introduce immediate extensions of
  valued modules and develop the corresponding Kaplansky theory as in 
  \cite{kaplansky}.
  \begin{defn}
  Given an extension $M\subseteq N$ of 
  $R$-modules, $N$ is said to be transcendental over $M$, if there is $x \in N$, such 
  that for all non zero $r$, $x.r \notin M$. Otherwise we say that $N$ is affine over
  $N$. In other words, every element of the $R$-module $N/M$ is a torsion. 
  \end{defn}


\subsubsection{Immediate Extensions}
Let $\Omega$ be a limit ordinal and $\{a_\gamma\}$ $(\gamma \in \Omega)$ a sequence from $M$, 
we say that $\{a_\gamma\}$ is a pseudo-Cauchy sequence if for all $\gamma >\delta>\rho$, sufficiently big, we have 
$$ v(a_\gamma - a_\delta)>v(a_\delta-a_\rho) .$$
We say that $\{a_\gamma\}$ pseudo-converge to some $z$ (possibly in an extension of $M$), and denote
by $\{a_\gamma \}\leadsto z$, if $v(z-a_\gamma)=v(a_{\gamma + 1} - a_{\gamma})$ eventually.
\begin{lemma} Let $\{a_\gamma\}$ be a pc-sequence and $r \in R$, be non zero, then $a_{\gamma +1} - a_\gamma$ is eventually regular for $r$ in any valued module. 
	\begin{proof}
		The difference $a_{\gamma + 1} - a_\gamma$ is eventually strictly increasing with $\gamma$ and hence its valuation is not a potential jump of $r$  (in particular not a jump value of $r$) since potential jumps of $r$ is a finite set. 
	\end{proof}
	\begin{corr}
		Set $z_\gamma:=a_\gamma.r - b$ for all $\gamma$. Then $\{z_\gamma \}$ is a pc sequence and 
		if $\{a_\gamma\} \leadsto a$ then $\{z_\gamma\} \leadsto a.r -b$ ($r\neq 0$).
	\end{corr}
	\begin{proof}
		We have $z_{\gamma}-(a.r-b)=(a_{\gamma}-a).r$. Eventually, $v(a_{\gamma}-a)$ is strictly increasing and by the above lemma $a_{\gamma}-a$ is regular  for $r$ .  Hence $v((a_{\gamma}-a).r)=v(a_{\gamma}-a)\cdot r$ is strictly increasing with $\gamma$ eventually. 
	\end{proof}
\end{lemma}
\begin{fact}The following is standard and is a direct consequence of the axioms of valued modules. 
  \begin{enumerate}
  \item Suppose $\{a_\gamma\} \leadsto z \in M$. Set $\delta_\gamma:=v(a_{\gamma+1}-a_\gamma)$. Then the set of all pseudo-limits of $\{a_\gamma \}$ in $M$, is 
  $$z + \{y \in M \tq v(y)>\delta_\gamma, \, \text{for all}\, \gamma \}.$$ 

\item Every pseudo-Cauchy sequence has a limit in an elementary extension of $M$.

\item Note that, $\{a_\gamma\} \leadsto 0$ if and only if 
  $v(a_\gamma)$ is eventually strictly increasing. Otherwise, if $z\neq 0$
  (in an possibly elementary extension) 
  a pseudo limit of the $\{a_\gamma\}$, then, 
  $v(a_\gamma)=v(z)<v(a_{\gamma} -z )$ eventually. 

  In conclusion the sequence $\{v(a_\gamma)\}$ is either eventually constant or is strictly increasing.
  \end{enumerate}
\end{fact}

\begin{defn}
  We say that $\{a_\gamma\}$ is of affine type if  for some $b\in M$ and for some non zero $r\in R$,
  $\{a_\gamma.r\} \leadsto b$. Otherwise we say that $\{a_\gamma\}$ is of transcendental type.
\end{defn}
\begin{remark}
  The sequence $\{a_\gamma\}$ is of transcendental type if and only if for all non zero $r
  $ and $b\in M$ the value $v(a_\gamma.r -b)$ is eventually constant and is of affine type if for some $b\in M$ and for some non zero $r$, $v(a_\gamma.r -b)$ is eventually strictly increasing with $\gamma$.

\end{remark}
  \begin{defn} Let $\{a_\gamma\}$ be of affine type. Then any non zero polynomial $r$ of minimal
  degree such that $\{a_\gamma.r\}$ has a limit in $M$ is be called a minimal polynomial
  of $\{a_\gamma\}$. 
  \end{defn}

\begin{defn}
Let $(M\subseteq N,v)$ be an extension of valued modules. If for every $x\in N$ there is
$y\in M$ such that $v(x-y)>v(y)$ then $(N,v)$ is called an immediate extension of $(M,v)$ We also say that the extension $(M\subseteq N,v)$ is immediate. $(M,v)$ is said to be maximal (respectively affinely maximal) if it has no proper immediate (respectively no
proper affine immediate) extension.
\end{defn}

\begin{remark} If the extension $(M\subseteq N,v)$ is immediate,
	then $v(M)=v(N)$ and  for every $\gamma$, 
	$N_{\gamma}=M_{\gamma}$, where $N_{\gamma}$ and $M_{\gamma}$ are respectively the quotients, in $N$ and in $M$, of the closed ball of radius $\gamma$ centered at $0$, by the open ball of radius $\gamma$ centered at $0$.
\end{remark}

The following is classical.
\begin{proposition}
	If $(M\subsetneq N,v)$ is immediate then every $b\in N\setminus M$ is a psedo-limit of some p.c. sequence from $M$,
	without limit in $M$.
\end{proposition}
\begin{proof}
	See proposition 4.3.13 in \cite{gonenc}.
\end{proof}

\begin{corollary}
If every pseudo-Cauchy sequence from $M$, has a limit in $M$, then $(M,v)$ is maximal.
\end{corollary}

By a valued module embedding (resp. isomorphism, automorphism etc.) we mean an $L$-embedding. That is, a map $(f,f_v)$ where $f$ is an $L_{Mod}(R)$-embedding between module sorts and 
$f_v$ is an $R$-chain embedding between value sets such that $v\circ f_v=w\circ f$.

The following propositions will establish a converse for the above corollary.
\begin{prop}\label{extimmtrans}
  Let $(M,v)$ be a valued modules and $\{a_\gamma\}$ a pc-sequence of transcendental type
  without limit in $M$. Then, $M(x):=M\oplus x.R$, where $x.R$ is the free $R$-module with 1-generator, admits a unique valued module structure (up to a valued module isomorphism fixing pointwise $M$) such that
  the extension $M \subset M\oplus x.R$ be immediate and $x$ be a pseudo-limit of 
  $\{a_\gamma \}$.
\end{prop}
\begin{proof} We define a the map $\tilde{v}:M\oplus x.R \to v(M)$, by setting
  $$\tilde{v}(m+x.r):=v(m + a_\gamma.r)$$ where $\gamma$ is big enough such that
  $v(m+a_\gamma.r)=v(m+a_\rho.r)$ for all $\rho\geq \gamma$. It is straight forward to
  check that $\tilde{v}$ satisfies the ultrametric inequality and that this yields a valued abelian group structure on $M(x)$. Now let $\mathbf{m}$ be
  a monomial of $R$. We have $\tilde{v}((m+x.r).\mathbf{m})=
  v(m.\mathbf{m} + a_{\gamma}.r\mathbf{m})$ eventually. Since monomials has no potential jump, 
   $v(m.\mathbf{m} + a_{\gamma}.r\mathbf{m})=\tilde{v}(m+a_\gamma.r)\cdot \mathbf{m}$
   eventually. It follows that $\tilde{v}((m+x.r).\mathbf{m})
   =\tilde{v}(m+x.r)\cdot \mathbf{m}$. Hence $(M(x), \tilde{v})$ is a valued module. 

  We will see now that $M':=M\oplus x.R$ is an immediate extension of $M$. 
  Fix $\gamma>\rho$ big enough such that such that
  $v(a_\rho-a_\gamma)>v(a_\gamma)=v(a_\rho)=v(x)$ (this is possible since $\{a_\gamma\}$ has no 
  limit in $M$, in particular $\{a_\gamma\} \not \leadsto 0$)
  and $\tilde{v}(x-a_\rho)=v(a_\gamma - a_\rho)$. Hence
  $\tilde{v}(x-a_\rho)=v(a_{\rho +1} - a_\rho)$ eventually for all $\rho$. This
  shows that $\{a_\gamma\} \leadsto x$. Now for $\rho$ big enough
  $\tilde{v}(m+x.r -m-a_\rho.r)=\tilde{v}((x-a_\rho).r)\geq \tilde{v}(x-a_\rho)\cdot r
  >v(a_\rho)\cdot r\geq \tilde{v}(m+x.r)$. Hence $M\oplus x.R$ is an immediate extension of $M$.
  
  Finally, uniqueness follows by the construction.
\end{proof}

\begin{corr}
  Let $(M \subseteq N,v)$ be an immediate extension and $y \in N\setminus M$ such that $y.r \in M$ for some non zero $r$. Then, $y$ is a limit of a pc-sequence of affine type without limit in $M$. 
\end{corr}

\begin{prop}\label{ext:imm}
  Let $(M,v)$ be a valued module and $\{a_\gamma\}$ be a pc-sequence of affine type
  without limit in $M$. Let $q$ be a minimal polynomial of $\{a_{\gamma}\}$ and $a\in M$ such that $\{a_\gamma.q\} \leadsto a$. Then $\deg(q)\geq 1$ and $(M\oplus x.R)/(x.q-a).R$ admits a unique valued module structure (up to a valued module isomorphism fixing pointwise $M$) such that
  the extension $M \subseteq (M\oplus x.R)/(x.q-a).R$ be immediate and $x$ be a pseudo-limit of $\{a_\gamma \}$.
\end{prop}
  \begin{proof}
Suppose $\deg(q)=0$, that is $q\in K^{\times}$. Then $v((a_\gamma.q -a.).q^{-1})=v(a_\gamma.q -a.)\cdot q^{-1}=v(a_\gamma - a.q^{-1})$ is strictly increasing with $\gamma$. That is $(a_{\gamma})_{\gamma}  \leadsto a.q^{-1} \in M$. This is a contradiction since we assume that $\{a_{\gamma}\}$ has no limit in $M$. 

Denote by $Z$ the submodule $(x.q-a).R$ of $M(x):=M\oplus x.R$. In $M(x)/Z$ every class has a unique representative of the 
for $x.s+m$ where $m\in M$ and $s$ of degree $<\deg(r)$. We set 
$\tilde{v}(x.s + m+ Z)$ as the ultimate value of the $v(a_{\gamma}.s +m)$. It is straightforward to check that  
$(M(x)/Z,\tilde{v})$ is a valued abelian group and 
$\tilde{v}(y.a)=\tilde{v}(y)\cdot a$ for all $a\in K$ and $y\in M(x)/Z$. 

We will now show that $\tilde{v}(y.t)=\tilde{v}(y)\cdot t$ for all  $y\in M(x)/Z$ and hence  establish that $(M(x)/Z,\tilde{v})$ is a valued module. 
Write  $y=x.s - m +Z$ where $m\in M$ and $\deg(s)<\deg(q)$. The case where $\deg(st)<\deg(q)$ follows immediately since $v(a_{\gamma}.st - m.t)$ is then ultimately constant. We may suppose  $\deg(st)=\deg(q)$.

 Since $v(a_{\gamma}.st-m.t) = v(a_{\gamma}.s-m)\cdot t$ and that the sequence
 $v(a_\gamma.s - m)$ is ultimately constant,
 the sequence 
$(v(a_\gamma.st - m.t))$ is ultimately constant.  Write  $st=q\mu + q_0$ with $\deg(q_0) < \deg(q)$ and $\mu \in  K^{\times}$.  Hence the canonical representative of $x.st -m.t +Z$
is $a.q_0 - m.t + a.\mu$.   
We have $a_{\gamma} .st-m.t = a_{\gamma} .q\mu-a.\mu+a.\mu+a_{\gamma}.q_0 -m.t$.
Note that $\{a_{\gamma} .q\mu-a.\mu \leadsto\} 0$ since $\{a_{\gamma}\} \leadsto a$. On the other hand, since $v(a_{\gamma} .st-m.t)$ is eventually constant,
$$v(a_{\gamma} .st-m.t)=v(a.\mu+a_{\gamma}.q_0 -m.t)$$ eventually. Now by definition, for all $\gamma$ big enough $v(a_{\gamma} .st-m.t)=\tilde{v}(x.s -m + Z)\cdot t$. 

The fact that the extension is immediate and unique up to an isomorphism of valued modules follow exactly as above. 
  \end{proof}
  
  Hence we have proved:

\begin{theorem}\label{thm:affmaxpc} A valued module $(M,v)$ is affinely maximal (respectively maximal), if and only if every pc-sequence of affine type (respectively every pc-sequence) has a limit inside $(M,v)$.
\end{theorem}

\begin{corr}\label{corr:appfields} Let $(F,v)$ be an henselian valued field of positive characteristic $p>0$. Let $R:=F[t;x\mapsto x^p]$,  identified with the ring
  of additive polynomials over $F$. Then the following assertions are equivalent:
  \begin{enumerate}[1.]
  \item $(F,v)$ is affinely maximal as a valued $R$-module.
  \item $(F,v)$ is algebraically maximal (as a valued field).
  \item For every one-variable usual polynomial $G$ over $F$, the set 
  $\{v(G(a)) \tq a \in F\}$ has a maximum in $vF$.
  \item For every additive polynomial $Q$ over $F$ and $b\in F$ the set 
  $\{v(Q(a)-b) \tq a \in F\}$ has a maximum in $vF$.
  \end{enumerate}
\end{corr}
\begin{proof}
  Suppose $(F,v)$ is not algebraically maximal. Take a proper immediate algebraic
  extension $(F\subseteq L,v)$. Then this extension is immediate in the sens of valued modules
  and $L$ is an affine extension of $F$ since every 
  polynomial divides (in the ring $F[X]$) an additive polynomial by \cite{ore}, Chapter 3, Theorem 1. Hence
  we get $1. \to 2.$ The assertion $2 \to 3$ is a folklore and the assertion
  $3 \to 4$ is trivial. 

  Let us show that $4 \to 1$. Suppose $(F,v)$ is not affinely maximal.  Then by Theorem \ref{thm:affmaxpc},
  there is a pc-sequence $\{a_\gamma\}$ of affine type, without limit in $F$,
  there is $b\in F$, such that $\{a_\gamma.q\}\leadsto b$, where $q$ is a minimal polynomial of $\{a_\gamma\}$. Hence $\{a_\gamma\} \leadsto \tilde{a} \in F^{alg}$
  with $\tilde{a}.q=b$. Let $a\in F$. Since $(F,v)$ is henselian, the extension of $v$ to the algebraic
  closure is unique and we have 
  $v(a.q-b)=\sum v(a-\tilde{a}_i) $ where $\tilde{a}_0=\tilde{a}$ and $\tilde{a_i}$ are the 
  $F$-conjugates of $\tilde{a}$. Since $F[\tilde{a}]$ is an immediate extension of $F$, there is
  $a_1 \in F$ such that $v(a_1-\tilde{a})>v(a-\tilde{a})$. Since $(F,v)$ is henselian,
  all the $a-\tilde{a}_i$  have the same valuation and $v(a_1.q-b)=\sum v(a_1-\tilde{a}_i)>v(a.q-b)$. That is $\{v(a.q-b) \tq a \in F\}$ has no maximum.
\end{proof}

\subsection{Residually divisible valued modules}
Residually divisible valued modules will be shown to be the analogs of the Kaplansky fields.
\begin{defn}\label{def:resdivisible}
  A valued module $(M,v)$ is called residually divisible if,
  \begin{enumerate}
  \item for all non zero $z\in M$ there exists $y$ such that 
  $v(y.t-z)>v(z)$,
  \item for all $r\in R$ with a unique potential jump $\gamma$, 
  for all $z$ of valuation $\gamma.r$, there exists $y$ of valuation
  $\gamma$ such that $v(z-y.r)>v((z)$.
  \end{enumerate}
\end{defn}
\begin{remark} By the above point (1),
  by induction and by axioms of valued modules, we have for every non zero monomial $\mathbf{m}$, for all non zero $z$ there exits some $y$ such that $v(y.\mathbf{m} -z)>v(z)$; in particular, $v(y.\mathbf{m})=v(y)\cdot \mathbf{m}=v(z)$. It follows that for every non zero $r$, the map $\gamma \mapsto \gamma\cdot r$ is a bijection of the value set. 
\end{remark}

The above two axioms are sufficient to get: 
\begin{lemma}\label{lem:regresdivisible} $(M,v)$ is residually divisible if and only if for all non zero $r \in R$ and non zero $z\in M$,
  there exists a regular $y$ for $r$ such that $v(y.r-z)>v(z)$. 
\end{lemma} 
\begin{proof} 
  Suppose $(M,v)$ is residually divisible. Let $r\in R$ and $z\in M$ be both non zero. 
  By the above remark, for some $\gamma$, $\gamma \cdot r=v(z)$. 
  Let $r'$ be the sum of the monomials
  $\mathbf{m}_i$ of $r$ such that $\gamma \cdot \mathbf{m}_i = \gamma\cdot r$. Hence 
  $\gamma \cdot r' =\gamma \cdot r$ and $r'$ is either a monomial or $\gamma$ is the unique 
  potential jump of $r'$. In both case there exists $y$, of valuation $\gamma$ such that 
  $v(y.r'-z)>v(z)$. Since $v(y.(r-r'))\geq v(y)\cdot (r-r')>v(y)\cdot r$, we get 
  $v(y.r-z)>v(z)$. It is easy to see that since $y$ is regular for $r'$ it is for $r$.

  The converse is obvious. 
\end{proof}
\begin{remark}\label{rem:resdiv<n} The above proof shows that if $(M,v)$ satisfies the 
  condition (1) of Definition \ref{def:resdivisible} and 
  the condition (2) for every polynomial of degree $\leq n$, then
  for all non zero $r$ of degree $\leq n$ and for all 
  non zero $z$ there exists a regular
  $y$ such that $v(y.r-z)>v(z)$.
\end{remark}

\begin{remark}
	Any immediate extension of a residually divisible module is residually divisible and a divisible module is residually divisible.
\end{remark}

\begin{defn} 
  A valued field of positive characteristic $p>0$, is said
  to be Kaplansky, if its value group is 
  $p$-divisible and it residue is field is $p$-closed, i.e. the 
  residue field $k$, is a divisible $k[t;x\mapsto x^p]$-module.
\end{defn}

\begin{lemma} A valued field of positive characteristic $p>0$ is Kaplansky if and only if it is residually divisible as a $K[t,x\mapsto x^p]$-module. 
\end{lemma}
\begin{proof}
Let $(K,v)$ be a Kaplansky field. Since $vK^{\times}$ is $p$-divisible, the map $\gamma \mapsto \gamma \cdot r$ is a bijection of the $R$-chain structure of $vK$. Let $r \in R$ be non zero. Let $z\in K$ be non zero, and $\gamma$ such that $\gamma\cdot r =v(z)$. Denote again by $r$ the additive polynomial associate to $r$.  Choose $u$ be of valuation $\gamma$, regular for $r$; we can always find a regular $u$ for $r$ by Remark \ref{regularelementsexits} since $k$ is infinite. Consider $f(x):=r(xu)/z -1$. It follows that $f$ has coefficients over the valuation ring $\mathcal{O}$ of $K$, and the induced polynomial $\bar{f}$ on $k$ is non constant (because $u$ is regular for $r$). Since $k$ is $p$-closed, there is some $w \in \mathcal{O}^{\times}$ such that $\bar{f}(\bar{w})=0$. It follows that $v((r(wu)-z)>v(z)$ and $wu$ is regular for $r$, since $r(wu)=v(z)=\gamma \cdot r$ and 
$v(uw)=\gamma$.

The converse is obvious.
\end{proof}
The proof of the following remark is exactly as above, by replacing $x\mapsto x^p$ by $x\mapsto x^{\sigma}$, in a valued difference field.
\begin{remark}\label{rem:resdivisible}
	Let $(K, \sigma, v)$ be an valued difference field, such that $\sigma$ is an automorphism and, the residue field $k$, is $\sigma$-linearly closed, that is every linear $\sigma$-polynomial over $k$, is surjective as a map $k\to k$. 
	Then the valued module $(K,v)$ is residually divisible as a $K[t;\sigma]$-valued module.
\end{remark}

\begin{theorem}\label{thm:aff_max} A valued module $(M,v)$ is residually
	divisible and affinely maximal  if and only if for all non zero $z\in M$ and all non zero $r$ 
  there exists a regular $y\in M$ for $r\in R$ such that $y.r=z$. In particular,
  $(M,v)$ is divisible whenever it is residually divisible and 
  affinely maximal. 
\end{theorem}
\begin{corr}\label{corr:aff_max} Let $(M,v)$ be residually divisible and affinely maximal. Then $x$ regular for $r$ (respectively regular), if and only if
	$\max_{a \in M}\{v(x-a) \tq a.r=0\}=v(x)$   (respectively 
	$\max\{v(x-a) \tq a \in M_{tor}\}$ exists and equal to $v(x)$). In particular,  any irregular $x$ for $r$, can be written as $x=a_0 + \epsilon$ where $v(\epsilon)>v(x)=v(a_0)$, $\epsilon$ is regular for $r$ and 
	$a_0.r=0$.
\end{corr}
\begin{proof}
	Let $x$ be irregular for $r$ and take $z$ regular for $r$ such that $x.r=z.r$. Then
	$(x-z).r=0$ and since $x$ is irregular and $z$ is regular for $r$, we $v(x)\cdot r< v(x.r)=v(z)\cdot r$. Hence $v(z)=v(x-(x-z))>v(x)$. Conversely, let $y$ be regular for $r$, then for any $a$ such that $a.r=0$, $v(y-a)\cdot r\leq v((y-a).r)=v(y.r)=v(y)\cdot r$, hence 
	$v(y)\geq v(y-a)$. 
\end{proof}

{\it Proof of the Theorem \ref{thm:aff_max}}.
  Suppose $(M,v)$ is affinely maximal and residually divisible. Let $z \in M$  and $r \in R$ be both non zeros. Suppose there is no $y$ regular for $r$ such that $z=y.r$. By Lemma $\ref{lem:regresdivisible}$ there exists $y_0$ such regular for $r$ such that $v(y_0.r-z)>v(z)$. 
  Since $y_0.r \neq z$, there exists $y_1$ regular for $r$, such that 
  $v(y_1.r -(y_0.r-z))=v((y_0-y_1).r -z)>v(y_0.r-z)>v(z)$. Note that, since $y_1$ is regular for $r$, and $\cdot r$ is strictly increasing, 
  $v(y_1)>v(y_0)$. Hence $y_0-y_1$ is also regular for $r$. Define
  by induction $y_n$ for $n \in \mathbb{N}$, such that 
  $z_n:=y_0 - (\sum_{i=1}^n y_i)$ be regular for $r$ with $v(z_n.r-z)>v(z_{n-1}.r-z)$ for all $n>0$. Then ${z_n}$ is a pc-sequence 
  of affine type with a limit in $M$ by hypothesis. Denote this limit by $z_{\omega}$. It is easy to check that $z_{\omega}$ is regular for $r$ and 
  $v(z_{\omega}.r -z)>v(z_n.r-z)$ for all $n$. But then we can continue to build
  $z_\lambda$ for any ordinal $\lambda$ in the same way. This is a contradiction considering the cardinality of $M$.

  For the converse, 
  let $\{a_\gamma\}$   be a pc-sequence of affine type from $M$ with a minimal polynomial 
  $r$ such that $\{a_\gamma.r\} \leadsto b \in M$. 
  We can write $a_\gamma=c_\gamma + d_\gamma$, ($c_\gamma, d_\gamma \in M$), where $v(c_\gamma) < v(d_\gamma)$ if $c_\gamma \neq 0$ and $d_\gamma$ regular for $r$ with $d_\gamma.r=a_\gamma.r$ and $c_\gamma.r=0$. By Lemma \ref{rvregular}, either $a_\gamma$ is eventually regular for $r$, either $a_\gamma$ is eventually irregular for $r$ (since $v(a_\gamma-a_\delta)>v(a_\gamma)$ for all $\delta>\gamma$ big enough).
   
   {\it Case 1}. $a_\gamma$ is eventually regular for $r$ and $b=0$:
   We have $\{a_\gamma.r\} \leadsto 0$. Since $\cdot r$ is strictly increasing and $v(a_\gamma.r)=v(\gamma)\cdot r$ eventually, we should have
   $\{a_\gamma\} \leadsto 0$. 
  
   {\it Case 2}. $a_\gamma$ is eventually irreguar for $r$ and $b=0$: Write $a_{\gamma+1}-a_\gamma=(c_{\gamma+1}-c_\gamma)-
   (d_{\gamma+1}-d_{\gamma})$. Since the $d_{\gamma}$ are regular for $r$, $d_{\gamma} \leadsto 0$ as in the {\it Case 1}. That is $v(d_{\gamma+1})>v(d_\gamma)$. Since   $a_{\gamma+1}-a_\gamma$ is eventually regular for $r$, we should have $v(c_{\gamma+1} -c_\gamma)\geq v(d_{\gamma})$ eventually.  Hence $c_{\gamma+1}-c_\gamma \leadsto 0$. This can only be possible if $c_{\gamma+1}=c_{\gamma}$ eventually, since $v(c_{\gamma+1}-c_{\gamma}) \in \jumpm{M}{r}$ which is finite. Hence, for some $c_0$, such that $c_0.r=0$,  eventually  $a_\gamma=c_0 + d_\gamma$, It follows that $a_\gamma \leadsto c_0$.
  
   {\it Case 3.}  $b \neq 0$: Let $b' \in M$ such that $b'.r=b$. It follows that $(a_\gamma -b').r \leadsto 0$. Working with the pc-sequence 
  $a_\gamma -b'$ we are in one of the above situations.  \qed

\section{Embedding Theorems and Quantifier elimination}
We will show a quantifier elimination result via the back and forth method. First we expose some preliminary observations.

Let $(M,v)$ be a residually divisible affinely maximal valued module. Note that $vM$ is in particular a $D$-chain. Let $Q$ be the set of irreducible elements of $R$ having a unique potential jump in the full closure of 
$vM$. We denote by $\jump{r}$, the set of jump values 
of $r$ in this full closure (recall the definition of full closure from \ref{def:fullchains}). 
 
Given $r\in R$ and $\gamma \in \jump{r}$,  we let $r_\gamma$ be the polynomial 
$$r_\gamma=\sum_i \mathbf{m}_i$$ where the $\mathbf{m}_i$ are exactly the monomials of $r$ such that $\gamma \cdot \mathbf{m}_i =\gamma \cdot r$. 
Note that,  $r_\gamma$ has at least two monomials. 

Now for $\gamma \in \jump{r}$, we denote by $\mathfrak{a}_{\gamma}$ the set of zeros of $r_\gamma$ in $M$. Let $A:=\{x\in M \tq x.r=0\}$. For  $\gamma \in \jump{r}$ we let 
$A_{\gamma}$ to be  a  -fixed- vector space complement of the 
$\Fix(\vfi)$-vector space $\{x \in A \tq v(x)> \gamma\}$ inside 
$\{x \in A \tq v(x)\geq \gamma\}$. In particular for any distinct $x,y \in A_{\gamma}$, $v(x-y)=\gamma$.

\begin{lemma}\label{jumpsrgamma}
	For every  $\infty \neq\gamma \in \jump{r}$, there exists a unique isomorphism of $\Fix(\vfi)$-vector spaces, $\xi_{\gamma}: \mathfrak{a}_{\gamma} \to A_{\gamma}$, such that
	$$v(\xi_{\gamma}(x)-x)>\gamma$$ 
	for all $x\neq 0$. Consequently, for all non zero $a \in A$ of valuation $\gamma$ there is a unique $a_\gamma$, root of $r_\gamma$, such that $v(a-a_\gamma)>\gamma$.

\end{lemma}
\begin{proof} For the uniqueness, note that  for every non zero $x \in \mathfrak{a}_{\gamma}$ and for two such map $\xi_{\gamma}, \xi'_{\gamma}$, we have
	$v(\xi_{\gamma}(x)-\xi'_\gamma(x))>\gamma$. But this means 	$\xi_{\gamma}(x)=\xi'_\gamma(x)$ since the unique element of $A_{\gamma}$ of valuation $>\gamma$ is $0$.
	
	Let $x\in \mathfrak{a}_{\gamma}$ be non zero. 	Then $x$ is irregular for $r$ and hence by Corollary \ref{corr:aff_max}, there is $y$, a zero of $r$ such that
	$v(x-y)>v(y)$. Note that such a $y$ is unique since if $z$ is another zero of $r$ of valuation $\gamma$ such that $v(x-z)>v(z)$ then $v(x-z)>0$ and hence
	$v(y-z)>\gamma$. In other words, $y=z$. We set $\xi_{\gamma}(x):=y$. Converse follows similarly since any non zero element of $A_{\gamma}$ is irregular for $r_{\gamma}$.

\end{proof}
\begin{remark}
 The above proof shows that 
$\mathfrak{a}_{\gamma}$ is trivial if and only if $A_{\gamma}$ is trivial.
\end{remark}

\begin{corollary}\label{cor:jumpsrgamma}
	$$\vert A \vert = \prod_{\gamma \in \jump{r}} \vert A_{\gamma} \vert = \prod_{\gamma \in \jump{r}} \vert \mathfrak{a}_{\gamma} \vert =\prod_{\gamma \in \jumpm{M}{r}} \vert \mathfrak{a}_{\gamma}\vert$$
	where $|\cdot| \in \mathbb{N}\cup \{\infty\}$.
	
\end{corollary}

Together with Corollary \ref{corr:aff_max}, we have also the following consequence that we will use very frequently in the rest of the paper. Hence we prefer not to referring every time. 

\begin{corollary}\label{corr:reg}
Let $x$ be non zero irregular element for some $r$, of valuation $\gamma$,  then there is a unique $a\in \bigoplus_{\delta\geq \gamma} \mathfrak{a}_{\delta}$, such that  $x-a$ is regular for $r$. 
\end{corollary}
\begin{proof}
By Corollary \ref{corr:aff_max}, there is a root $a_0$ of $r$ such that $v(x-a_0)>\gamma$. By above Lemma, there is a unique $a_\gamma$ root of $r_\gamma$ such that $v(a_0-a_\gamma)>\gamma$. Hence $\delta:=v(x-a_\gamma)>\gamma$. Now if $x-a_\gamma$ is irregular, then repeating the same argument there is a unique $a_\delta$, root of some $r_\delta$, such that $x-a_\gamma - a_\delta$ is of valuation $>\delta$. Since $\jumpm{vM}{r}$ is finite, either we find $a$ such that $v(x-a)>\max \jump{r}$ and hence $x-a$ regular for $r$, either we already obtain $a$ such that $x-a$ is regular for $r$  in less than $\jump{r}$ step. 
\end{proof}

Recall that $r \in R$ is said to be separable if $t$ does not divides $r$ and an irreducible $r$ is separable whenever $r\neq t$. Moreover every $r$ can
be written as $r=t^ns$ where $s$ separable. In particular 
$$\vert \{x \in M \tq x.r=0\} \vert =
\vert \{x\in M \tq  x.s=0\} \vert.$$

\begin{prop}
Let $(M,v)$ be a residually finite and affinely maximal valued module. Then the function $\eta_M:R\setminus\{0\} \to \mathbb{N}\cup \{\infty\}$, 
$r \mapsto |\{x\in M \tq x.r=0\}|$ is determined by its
restriction on $Q$.
\end{prop}
\begin{proof}
Let  $ d \in \mathbb{N}\setminus{0}$ and $R_{<d}$ be the set of non zero polynomials of degree $< d$. We proceed by induction. Let $r \in R$ be of degree $d$. We may suppose that $r$ is separable. If $d=1$ the $r$ is necessarily  irreducible and has at most 2 monomials. 

Suppose $d > 1$ and the restriction of
$\eta_M$ to  $R_{<d}$ is determined by the restriction of $\eta_M$ to $Q$. We write the irreducible decomposition of $r$ as
$r = r_1 \ldots r_n$. Let $A$ denotes the zero set of $r$ in $M$ and $A_i$ the zero set of $r_i$. Then by Lemma 2.17 in \cite{Onay2017}, we have $|A| = \prod_i |A_i|$. Then, if  $n > 1$, then we got the result by induction. 

We may hence suppose that $r \notin Q$ and is irreducible.  Since $r \notin
Q$, $|\jump{r}| >2$. By Corollary \ref{cor:jumpsrgamma},

$|A|=\prod_{\gamma \in \jump{r}}|\mathfrak{a}_\gamma|$. Let $\delta = \min \jump{r}$.
Then $r_{\delta}$ is not separable. Write $r_\delta=t^ks$ with $k>0$ and $s$ separable. Then $|\mathfrak{a}_{\delta}|=|\{x\in M \tq x.s=0\}|$. On the other hand for all $\gamma > \delta$, 
$r_\gamma$ is of degree $<d$. Hence we may apply the induction hypothesis.
\end{proof}

Let $(K,\vfi)$ be a difference field, and $R:=K[t; \vfi]$ as usual.

 Let $(F, v,\Delta,)$ be an affinely maximal residually divisible valued module.  We suppose that the valuation $v_K$ induced on $K$ is non-trivial. 
 
 Let $T_v$ be the complete  $L_V$-theory of $v(F)$ and  $\operatorname{Tor}_F$ 
the $L_{Mod}(R)$-theory consisting of true sentences 
of the form $|\eta_F (r)| = n$ ($n \in \mathbb{N}\cup\{\infty\},
r \in Q$). Let $L:=L_V\cup L_{Mod}(R)$. Note that by Theorem 1.2 \cite{Onay2017}, this the complete  $L_{Mod}(R)$-theory of $F$. 

We denote by 
$T_{M}$ be the $L$-theory of affinely maximal residually divisible $R$-modules. Set
$$T := T_M \cup T_v \cup \operatorname{Tor}_F.$$ 

\begin{theorem}\label{thm:main}
	$T$ eliminates quantifiers of the module sort and  is complete.
\end{theorem}

 The essential content of the proof is given in the  Proposition \ref{back-forth}. First we have two more preliminary lemmas.

Let $(M,v)$ and $(N,w)$ be residually divisible affinely maximal valued modules.

Suppose $A\subseteq M$ be a submodule with an $L_{Mod}(R)$-embedding $f:A \to N$ with $f(A)=B$, and $f_v:v(A)\to w(B)$ be an $L_V$-embedding such that $f_v\circ v(a) = w\circ f(a)$. Let $x\in M\setminus \{0\}$ be affine over $A$ with $r_0$ a minimal polynomial of $x$ over $A$. Let $y\in N$ be such that $y.r_0=f(x)$.

The proofs of following two lemmas are straightforward and detailed in Lemmas 4.6.6 and 4.6.7 in \cite{gonenc}.
\begin{lemma}\label{lem:prep1} 
Suppose for all non zero $r$ of degree $<\deg(r_0)$ and for all $z$,
 if $z.r \in A$ then $z \in A$ and that $g_v:v(A+x.R)\to w(N)$ be an $L_V$-embedding extending
$f_v$ and satisfying, for all $a\in A$ $$g_v(v(x-a))=w(y-f(a)).$$
Then the extension $g \supseteq f$ mapping $x\mapsto y$ yields that 
$(g,g_v)$ is a valued module embedding of $(A+x.R,v)$ into $(N,w)$.
\end{lemma} 
\begin{lemma}\label{lem:prep2}
  Suppose now $g_v:v(A+x.R) \to w(N)$ satisfying $w(y)=g_v(v(x))$ and suppose also 
  that $x.r_0\neq0$. Then 
  \begin{enumerate}
  \item if $a\in A$ and $x-a$ and $y-f(a)$ are both regular for $r_0$, then 
  $$g_v(v(x-a))=w(y-f(a)),$$
  \item if, $A$ contains all the roots of $r_0$ inside $M$, then the following are 
  equivalent
  \begin{enumerate}
  \item for all $a\in A$, $g_v(v(x-a))=w(y-f(a))$,
  \item for all $a\in A$, $x-a$ is regular for $r_0$ if and only if 
  $y-f(a)$ is regular for $r_0$
  \end{enumerate}
  \end{enumerate}
\end{lemma}


Recall from our first article \cite{Onay2017}, that if $C$ is a submodule of an divisible module $M$, then there is unique divisible closure of $C+M_{tor}$, which is the submodule 
$$(C+M_{tor})^{div}:=\{x\in M \tq x.r \in C \mesp \text{for some non zero}\mesp r\in R\}.$$ 

\begin{prop}\label{back-forth}
 Let $(C,\Delta)$ and $(C',\Delta')$ be  respectively substructures of  $(M,v)$ and $(N,w)$, both models of $T$,  $L$-isomorphic via 
 $\mathbf{f} = (f, f_v )$. Then $\mathbf{f}$ extends to an isomorphism 
 $$\hat{\mathbf{f}}:\mesp \left((M_{tor} + C)^{div}, \operatorname{dcl}_{L_V}(\Delta)\right) \xrightarrow{\simeq}
 \left(N_{tor} + C')^{div}, \operatorname{dcl}_{L_V}(\Delta')
 \right).$$
\end{prop}

\begin{proof} 
Note that $vC\subseteq \Delta$ and $vC' \subseteq \Delta'$. If $x.r \in C$, for some non zero $r$, then either $x$ is regular for $r$ and 
$v(x) \in \operatorname{dcl}_{L_V}(vC)$ since $v(x)\cdot r=v(x)\cdot t^ia$ for some integer $i$ and $a\in K$, 
either $v(x)\in \jumpm{M}{r}\subseteq \jump{r} \subseteq
\operatorname{dcl}_{L_V}(\emptyset)$. Hence $v(C+M_{tor})^{div}\subseteq 
\operatorname{dcl}_{L_V}(vC)\subseteq 
\operatorname{dcl}_{L_V}(\Delta)$. Since $vM \equiv wN$, and $vC \xrightarrow[f_v]{\simeq} wC'$, $f_v$ extends  uniquely  to an $L_V$-isomorphism $\operatorname{dcl}_{L_V}(\Delta) \to \operatorname{dcl}_{L_V}(\Delta')$. 

For the rest of the proof, we note this extension again by $f_v$ and we assume $\Delta=\operatorname{dcl}_{L_V}(\Delta)$ and  $\Delta'=\operatorname{dcl}_{L_V}(\Delta')$. Moreover, if
$f:A\to N$ is a partial isomorphism of $R$-modules, we say that $f$ respect $f_v$, if $w(f(a))=f_v(v(a))$ for all $a\in A$.

  We set $B:=(C+M_{tor})^{div}$. In the following lines, we will extend $f$ to $\hat{f}:=B \to N$ so that $\hat{f}$ respects $f_v$. Hence we will obtain a partial isomorphism extending $\mathbf{f}$ on $(B,\Delta)$.

Note that by Corollary \ref{cor:jumpsrgamma}, together with 
axioms $\text{Tor}_{F}$, $f_v$ induces an isomorphism of ordered sets 
between $\jumpm{M}{r}$ and $\jumpm{N}{r}$ for all non zero $r$.

We will achieve the proof by showing by induction over $n \in \mathbb{N}$ the following
assertion:\\[1.3mm]
For all submodule $A$, such that 
$B \supseteq A \supseteq C$, if  $f:A \to A' \subseteq N$ is an
isomorphism  respecting $f_v$, then it extends
respecting $f_v$, to a
submodule $\bar{A} \subseteq B$ such that, for all $q \in R$ of degree $\leq n$, for all $x\in M$,
  $$x.q \in \bar{A} \Rightarrow x \in \bar{A}.$$

Let $A$ given as above. The assertion is trivial for $n=0$ since $A$ is a $K$-vector space. 
Fix $n >0$. We split the proof in four steps.

{\bf Step  1.} The map $f$ extends respecting $f_v$, to a submodule $\hat{A}$, such
that for all $r$ of degree $n$,
all $x \in M$ of valuation $> \max \jumpm{M}{r}$ 
(note that each such $x$ is regular for $r$),
if $x.r$ is non zero and $x.r \in \hat{A}$ then $x \in \hat{A}$.

If $\jumpm{M}{r}$ is empty we convey that for all $x \in M$ such that $x.r \in A$,
$v(x) > \max \jumpm{M}{r}$.
Let $r$ be non zero and $x\in M\setminus A$ such that $x.r=a_0\in A$ with $v(x)>\max\jumpm{M}{r}$. Necessarily, $r$ is a minimal polynomial of $x$ over $A$. 
 Choose $y \in N$, regular for $r$ and $y.r = f (a_0)$. 
 Since
 $f_v (v(x.r)) = w(f (a)) = w(y.r)$ and $x, y$ are both regular for
 $r$, we have $f_v(v(x)) = w(y)$. In particular,
$ w(y) > \max \jumpm{N}{r}$. It follows that for all $a\in A$ 
 such that $v(a)\geq v(x)$, $x-a$ and $y-f(a)$ are both regular for $r$. Hence, by Lemma \ref{lem:prep2} (1), $f_v(v(x-a))=w(f(a)-y)$. Moreover, for $a\in A$, we have $v(a)<v(x)$ if and only if $w(f(a))<w(y)$. Hence $v(x-a)=v(a)$ and 
 $w(y-f(a))=w(f(a))$. That is, the hypothesis of Lemma \ref{lem:prep1} are satisfied. Hence we can extend $f$ to $A+x.R$ respecting $f_v$.

\begin{itemize}
\item[$\star$] By induction hypothesis there exists a submodule $A(x)_{<n} 
\supseteq A + x.R$ such that if
$x' \in M$ and $x'.s \in A(x)_{<n}$ for some non zero $s$ of degree $<n$ then 
$x'\in A(x)_{<n}$ and $f$
extends on it.
By applying the above process and by transfinite induction
we can construct a submodule $A_1 \supseteq A$ such that for all $r$ of degree $n$
and for all $x' \in M$, if $x'.r=a (\neq 0) \in A$ and if $v(x' ) > \max \jumpm{M}{r}$, then $x' \in A_1$ .
We can now construct for all integer $i>0$, $A_{i+1}$ from $A_i$ in the same way and extend $f$. 
Then $\hat{A}:=\cup_i A_i$ has the required propriety. 
\end{itemize}
\qed ({Step 1}).

For the rest of the proof we assume $A=\hat{A}$.

{\bf Step 2.} We will now extend $f$ (respecting $f_v$) to a submodule $A'\supseteq A$  which contains every $x \in M$ such that
$x.r \in A'$ for some $r \in R$ of degree $\leq n$, and such that $\jumpm{M}{r}$ is a singleton.

{\bf 2.1} We will first add to $A$ the roots of polynomials $r$, of degree $n$ such that $\jumpm{M}{r}$ is a singleton,
and extend $f$. Choose a such polynomial $r$ of degree $n$, and $x \in M\setminus A$ be such that $x.r =0$.
Let $y \in N$ a root of $r$. Then necessarily $f_v(v(x)) = w(y)$ by the discussion in the beginning of the proof. In addition,
for all $a \in A$, $v(x- a) = \min\{v(x), v(a)\}$. In fact, otherwise, $x-a$ is of valuation
$> \max \jumpm{M}{r}$ and $(x-a).r=-a.r \in A$. By preceding step $x-a \in A$,  hence $x \in A$: a contradiction.
This implies, for all $a \in A$, $f_v(v(x-a)) = w(y-f (a))$. 
By Lemma \ref{lem:prep1} we can extend $f$ by
sending $x$ to $y$, on $(A + x.R)$ respecting $f_v$. We then apply above process ($\star$).

{\bf 2.2} We will now add all $x \in M$ , $x.r \in A$ for some $r$ of degree $n$ such that $\jumpm{M}{r}$ is a singleton.

 Let $x$ and $r$ be as in hypothesis. Then, $v(x-a)$ is necessarily regular for $r$, for
all $a \in A$. Otherwise for some $a'$, $a'.r=0$ and $v(x -a -a' ) > v(x -a)$ and
$x - a - a'$ is of valuation $> \max \jumpm{M}{r}$, hence $x \in A$. In particular $x$ is regular for $r$. Then for any $y \in N$ such
that $y.r = f(x.r)$ and for any $a \in A$, $y-f(a)$ is regular for $r$. By Lemma \ref{lem:prep2}, $f$
extends on $(A + x.R)$ respecting $f_v$. We then apply above process ($\star$). 

\qed(Step 2.)

For the rest of the proof we assume $A=A'$. Note that at this stage, by \ref{lem:regresdivisible}, $A$ is in particular residually divisible with respect to polynomials
of degree $\leq n$.

{\bf Step 3.} 
 We will now add to $A$, the roots of arbitrary polynomials of degree $n$.  We will now extend $f$ to a submodule $A_0\supseteq A$  which contains every $x \in M$ such that
 $x.r=0$ for some $r \in R$ of degree $\leq n$.

Let $r$ be of degree $n$ and $x \in M \setminus A$ such that $x.r = 0$. Note that
the case $|\jumpm{M}{r}|=1$  is treated on the previous step.

We will show that $x$ is a pseudo-limit of a
pc-sequence $\{a_\gamma\}$ of affine type, without limit in $A$, with $r$ a minimal polynomial of this sequence. 

 This will imply that $A+x.R$
is an immediate extension of $A$ and  by Proposition \ref{ext:imm}, its isomorphism type
 is uniquely determined so
that we can extend $f$ on $A+x.R$ respecting $f_v$.

Let $\delta=v(x)$. Since $x$ is non zero, $\delta \in \jumpm{M}{r}$. Recall that the subpolynomial  $r_\delta$ is of degree 
$\leq n$ and has only one jump value. If $r=r_\delta$, that is $r$ has unique jump value. Then there is nothing to do.

Otherwise by {\bf Step 2},  the roots
of subpolynomial $r_\delta$ are in $A$. It follows  by Corollary \ref{corr:reg} that there is $a_0 \in A$ 
such that $x-a_0$ is regular for $r$. If $a_0.r=0$ then necessarily $x=a_0 \in A$ since $0$ is the unique regular root of $r$. Hence $a_0.r \neq 0$.

Since $A$ is residually divisible for polynomials of degree $n$, by Remark \ref{rem:resdiv<n}, there
exists $a_1' \in A$ regular for $r$ such that  $v(a_0.r+a_1'.r)>v(a_0.r)$. Hence $v(a_0)<v(x-a_0)=v(a_1')$.  Then either, already
$(x-a_0-a_1')$ is regular for $r$ and hence 
$v(x-a_0-a_1')>v(x-a_0)$ and we set $a_1=a_0+a_1'$; or there exists some $a\in A$ such that $v(a)=v(x-a_0)=v(a'_1)$ and $x -(a_0+a'_1+a)$ is regular for $r$, hence  and we set $a_1:=a_0+a'_1+a$. Note that
$v(a_1.r)=v((a_1-x).r)=v(a_1-x)\cdot r > v(a_0 -x)\cdot r=v((a_0-x).r)=v(a_0.r)$.
  
By induction  we can
construct $a_{k+1}$ for all $k\in \omega$, such that $x-a_{k+1}$ is regular for $r$ and 
$v(a_{k+1}.r)>v(a_k.r)$. Hence the sequence $\{a_k\}$ is a p.c. sequence in $A$, $\{a_k.r\} \leadsto 0$ 
and $\{a_k\} \leadsto x$.
If $\{a_k\} \leadsto a_{\omega} \in A$ then it is straightforward to check that
 $v(a_{\omega} - a_k)=v(a_{k+1}-a_k)$ and $v(a_{\omega}.r)>v(a_k.r)$ for all $k$. 
Note that  $a_{\omega}.r\neq 0$ since by the  uniqueness of $a_0$ we should have
$a_{\omega}=x$. Thus we can continue by choosing $a_{\omega+1}$,  $a_{\omega+2}, \dots$, so on;
until have been constructed a p.c. sequence without limit in $A$. 

We apply the process ($\star$). \qed(Step 3.)

For the rest of the proof we assume $A=A_0$.

{\bf Step 4.} It remains to add the set of $x \in M\setminus A$
such that $x.r \in A \setminus \{0\}$, for some $r$ of
degree $n$  with $|\jumpm{M}{r}|\geq 2$.

Let $(x,r)$ be a such pair. Remark that $\Delta:=\{v(a'-x) \tq a'\in A\}$ admits a maximum of the form $v(x-a_0)$ for some $a_0\in A$. Otherwise 
there is a p.c sequence $\{a_\gamma\} \leadsto x$, without limit in $A$, with $(x-a_{\gamma}).r \leadsto 0$. But then there is proper immediate extension of $A$, affine over $A$, containing a root of $r$, but this is absurd. 

Set $z:=x-a_0$ where $a_0 \in A$ is such that $v(x-a_0)=\max \Delta$. Notice that for every $a\in A$, such that $v(a)\geq v(z)$, 
$z-a$ is regular for $r$, since otherwise for some $a' \in A$, root of $r$, $z-a-a'$ would be of valuation
$>v(z-a)=v(z)$. Let $y\in N$ be regular for $r$ such that $y.r=f(z) \in B$. Then necessarily,  
$w(y)=\max\{w(y-b) \tq b\in B \}$. It follows that 
$f_v(v(z-a))=w(y-f(a))$ for all $a\in A$. Hence mapping 
$z\mapsto y$, extends $f$ to $A+z.R=A+x.R$, respecting $f_v$.
\qed (Proposition.)


We recall first the following before attacking the proof of Theorem \ref{thm:main}. 

As in our first article \cite{Onay2017} two submodules $(A \subseteq M,v)$ and $(A'\subseteq M,v)$, are said to be {\it in sum valuation independent}, if for all $a \in A$ and $a'\in A$, $v(a\pm a')=\min\{v(a),v(a')\}$. It follows that if $B,B'$ are submodules of $N$, such that $(B,w)$ and $(B,w')$ are in sum valuation independent and  ${\bf f}=(f,f_v):(A,v)\to (B,w)$ and ${\bf g}=(g,g_v):A'\to B'$ are valued module isomorphisms 
 then the map $(f\oplus g,\min(f_v,g_v))$ is a valued module isomorphism from $(A\oplus A', v) \to (B\oplus B',w)$.
\end{proof}

{\bf Proof of the theorem \ref{thm:main}.}

We suppose both $(M,v)$ and $(N,w)$ are $|K|^{+}$-saturated models of $T$. We will show that if $C:=(A,\Delta)$ is substructure of size $\leq |K|$, ${\bf f}=(f,f_v):(A,\Delta) \to (N,w)$ a partial isomorphism and $x\in (M,v)\setminus C$, 
we can extend ${\bf f}$ to substructure containing $x$ in its domain.

By the above proposition we may assume that $A$ contains $M_{tor}$ and is affinely maximal (in particular divisible).

Let $\gamma \in \Delta \setminus vA$ and $x\in M$ of valuation $\gamma$. Since $vA \supseteq \jumpm{M}{r}$ for all $r$, $x$ is regular for all $r\in R$. It follows that $\gamma \cdot r \notin v(A)$ for all $r$, since $A$ is divisible. It follows that  $A$ and $x.R$ are in sum valuation independent. Let $B=f(A)$. Let 
$\delta=f_v(\gamma) \in wN$ and $y \in N$ of valuation $\delta$, then necessarily $y$ is regular for all $r$, and $B$ and $y.R$ are in sum valuation independent. Hence $f$ extends on $A\oplus x.R$ by sending $x \mapsto y$, respecting $f_v$. 

Note that $v((A+x.R)^{div}) \subseteq \operatorname{dcl}_{L_V}(\Delta)$. In fact  $v((A+x.R)^{div})=v(A)\cup \gamma\cdot t^{\mathbb{Z}}K$.  By the proof of the above proposition, we can extend $f_v$ to $\operatorname{dcl}_{L_V}(\Delta)$ and then extend $f:
(A+x.R)^{div} \to (B+y.R)^{div}$ respecting $f_v$. 

It follows that for the rest of the proof we may assume that 
$v(A)=\Delta$, and $w(B)=\Delta'$.

Now let $x\in M\setminus A$ and $v(x)=\alpha \notin vA$.
 Note that the quantifier free $L_V(\Delta')$-type of $\alpha$ (this is the set of  quantifier free formulas $\phi(x,\delta')$, $\delta' \in \Delta'$, such that $M\models \phi(\alpha,f_v^{-1}(\delta')$) is finitely satisfiable in $\Delta'$. By saturation and by the fact that $f_v$ is an partial isomorphism preserving $L_V(\Delta)$-formulas, there exists $\beta \in w(N)$ satisfying 
the quantifier free type $L_V(\Delta')$-type of $\alpha$. Hence by sending $\alpha \mapsto \beta$, $f_v$ extends on
$vA\cup \{\gamma \cdot r \tq r \in R\}$. As above, for all $y\in N$ of valuation $\beta$, the submodules $B$ and $y.R$ are in sum valuation independent. For a such $y$, sending $x\mapsto y$ induces an extension of $f$ to $A\oplus x.R$ with image $B\oplus y.R$, respecting $f_v$.

Now by the above proposition, we can extends ${\bf f}$ to $((A+x.R)^{div},v)$.

We have now two cases which remains. Let $x\in M\setminus A$. 

1. Suppose there is $a \in A$, such that 
$$v(x-a)=\max\{v(x - a') | a' \in A\}.$$ Since $A$ contains $M_{tor}$, $x-a$ is regular for all  $r\in R$.  Moreover, 
$A$ and $(x - a).R$ are in sum valuation independent. But then by saturation, there exits $y \in N$ such that $w(y - f (a)) = \max\{w(y - b) | b \in B\}$. It follows that we can extend $f$ to  $A + (x - a).R$ by mapping $x \mapsto y$.
	
Now as is the proposition above we can extend $f$ to the divisible closure of $A+x.R$, respecting $f_v$.

2. We may now suppose that $A \subseteq M$ is immediate. 

Let $x\in M\setminus A$. Hence $x$ limit  of a pseudo-Cauchy sequence from $A$. A such sequence is necessarily of trancendant type since $A$ is affinely maximal. Hence the isomorphism type as a valued module of the extension $A \subset (A + x.R, v)$ is uniquely determined. Hence we can extend $f$ on $(A+x.R,v)$, respecting $f_v$. \qed

The following results are immediate consequences of Theorem
\ref{thm:main}.

\begin{intth}[A-K,E $ \equiv $, Theorem {\bf \ref{akeequiv}}]
	Let $ (F, v) $ and $ (G, w) $ be two non zero affinely maximal residually divisible modules, such that $F_{tor}$ and $G_{tor}$ are elementary equivalent as $L_{Mod}(R)$-structures and $ v (F) $ and $ v (G) $ are elementary equivalent in the language $ L_{V} $. Then $ (F, v) $ and $ (G, v) $ are elementary equivalent as $ L $-structures.
\end{intth}

\begin{theorem}[A-K,E $ \preceq $]\label{corr}
Let $ (F, v) $ and $ (G, w) $ be two non zero affinely maximal residually divisible modules, such that $F_{tor}$ elementarily embeds in $G_{tor}$  as $L_{Mod}(R)$-structures and $ v (F) $ elementarily embeds in $ v (G) $  in the language $ L_{V} $. Then $ (F, v) $ elementarily embeds in an elementary extension of $ (G, v) $ as $L$-structures.
\end{theorem}

\begin{corollary}\label{corr:ann_infinite}
	Suppose $\Fix(\vfi)$ is infinite then the two above results hold without the requirements on torsion submodules.
\end{corollary}  
\begin{proof}
	Since the annulators are $\Fix(\vfi)$-vector spaces, the condition on the size of annulators are automatically satisfied, that is we have $\eta_M(r)=\eta_N(r)=\infty$, for all $r\in R$.
\end{proof}

Now by Remark \ref{rem:acfa} and \ref{example-full-chains}(2) we obtain Theorem \ref{thm:acfa}.
\begin{theorem}[Theorem {\bf \ref{thm:acfa}}]
	Let $(K,v_K, \sigma)$ be a valued difference field such that $(K,\sigma) \models \operatorname{ACFA}_{p}$ and $v_K(x^{\sigma})>nv(x)$ for all $x$ such that $v_K(x)>0$. The theory of affinely maximal residually divisible $K[t;\sigma]$-modules together the theory of dense $R$-chains is complete, eliminates quantifiers in $L$, and the valued $R$-module $(K,v_K)$ is the prime model of this theory.  
\end{theorem} 
%
%
%
%
%
%
%
%
%
%

  \bibliographystyle{plain}
  \bibliography{onayg}

\end{document}